\documentclass[12pt]{amsart}
\usepackage{amsmath,amsfonts,amsthm, amscd, mathrsfs, amssymb, txfonts}
\usepackage{xy}
\usepackage{fullpage}

\theoremstyle{plain}

\numberwithin{equation}{section}
\newtheorem{theorem}{Theorem}[section]
\theoremstyle{definition}

\newtheorem{lemma}[theorem]{Lemma}

\newtheorem{proposition}[theorem]{Proposition}

\newtheorem{corollary}[theorem]{Corollary}

\newtheorem{conjecture}[theorem]{Conjecture}
\theoremstyle{definition}

\newcommand{\Supp}{\textrm{Supp}}
\begin{document}

\title{On a Conjecture of Hacon and McKernan in Dimension Three}\thanks{\noindent Shortly before this article was finished Viehweg and Zhang posted an article on the arXiv, arXiv:math.AG/0707.4287, with stronger results along the same lines of this paper, using different methods.}
\author{Adam T. Ringler}
\date{August 24, 2007}
\address{Department of Mathematics and Statistics, University of New Mexico, Albuquerque, New Mexico 87131}
\email{ringler@unm.edu} 

\maketitle

\begin{abstract}
We prove that there exists a universal constant $r_3$ such that if $X$ is a smooth projective threefold over $\mathbb{C}$ with non-negative Kodaira dimension, then the linear system $|r K_X|$ admits a fibration that is birational to the Iitaka fibration as soon as $r \geq r_3$ and sufficiently divisible.  This gives an affirmative answer to a conjecture of Hacon and McKernan \cite[Conjecture 1.7]{McKernan} in the case of threefolds. 
\end{abstract}

\section{Introduction}

A fundamental problem in complex projective geometry is to understand the linear system $|rK_X|$ on a smooth projective variety $X$ when the Kodaira dimension $\kappa(X)$ is non-negative.  Recently, Hacon and McKernan proved in \cite{McKernan} a fundamental result stating that there exists a universal constant $r_n$ such that if $X$ is a smooth projective variety of dimension $n$ over $\mathbb{C}$ and of general type, then the pluricanonical maps
$$
\phi_{|rK_X|} : X \dashrightarrow  \mathbb{P} (H^0 (X, \mathscr{O}_X (r K_X)))
$$
are birational as soon as $r \geq r_n$.  They also conjectured \cite[Conjecture 1.7]{McKernan} that a similar type of result should hold whenever $\kappa(X) \geq 0$.  

\begin{conjecture}\label{basic}
There exists a universal constant $r_n$ such that if $X$ is a smooth projective variety of dimension $n$ with $\kappa(X) \geq 0$, then the pluricanonical map
\begin{equation}\label{thing}
\phi_r=\phi_{|r K_X|} : X \dashrightarrow  \mathbb{P}(H^0(X, \mathscr{O}_X (rK_X)))
\end{equation}
is birational to the Iitaka fibration, whenever $r \geq r_n$ and sufficiently divisible.
\end{conjecture}

McKernan has pointed out that \cite[Conjecture 1.7]{McKernan} needs to include the hypothesis that $r$ is sufficiently divisible.  One would also like the divisibility criterion for $r_n$ to depend only on $n=\dim X$.  For example if the conjecture holds for $r_n >1$ but, we require that $r \geq r_n$ be divisible by $1,000,000$, or even worse the divisibility depends on invariants of $X$ which are not bounded, then we have not ``effectively computed" the constant $r_n$.  

In this paper we show that Conjecture \ref{basic} is true when $\dim X =3$.  We also show that it is possible to bound the divisibility requirement independently of $X$.  That is, we prove:

\begin{theorem}\label{dim3}
If $\dim X =3$ and $\kappa(X) \geq 0$, then there exists a constant $r_3$ such that
$$
\phi_r : X \dashrightarrow \mathbb{P}( H^0 (X, \mathscr{O}_X (r K_X)))
$$
is birational to the Iitaka fibration of $X$ as soon as $r \geq r_3$ and $r$ is sufficiently divisible.  Furthermore, the divisibility criterion can be bounded independent of $X$.
\end{theorem}

Based on the methods used in this paper, we are able to get a rather explicit bound for $r_3$ when $\kappa(X)=2$:

\begin{corollary}\label{kap2bound}
When $\kappa(X) =2$ and $r \geq 48$ and divisible by $12$, then $\phi_r : X \dashrightarrow \mathbb{P}(H^0(X, \mathscr{O}_X (rK_X)))$ is birational to the Iitaka fibration.
\end{corollary}

It is also worthwhile to point out that the constant $r_3$ can, in general, be bounded explicitly.  However, the bound when $\kappa(X)=1$ is difficult to compute and presumably far from optimal.

\begin{corollary}\label{explicit}
The bound appearing in Theorem \ref{dim3} and the divisibility criterion are effectively computable.
\end{corollary}

Unfortunately, the methods used to prove Theorem \ref{dim3} seem to be rather restricted to when $\dim X=3$.  Therefore, another approach will be needed to prove Conjecture \ref{basic} in the general case.  However, Pacienza has made some recent progress, with some additional hypotheses, in the higher dimensional case, see \cite{Pacienza}.  We should also point out that, at least in the case when $\dim X=3$, Koll{\'a}r has some results in this direction using other methods, see \cite[Section 7]{Kollar4}. 

We now give an outline of the proof of Theorem \ref{dim3}.  Based on previous work of Fujino, Hacon, McKernan, and Mori, see \cite[Theorem 1.1]{McKernan}, \cite[Corollary 6.2]{MoriFuj}, and \cite[10.4.3]{Mori}, it suffices to deal with the case where $\kappa(X)=2$.  After making a birational modification we can further assume that
$$
\phi_r : X \to S_r \subset \mathbb{P} (H^0 (X, \mathscr{O}_X (rK_X)))
$$
is an elliptic fiber space, where $\phi_r (X)=S_r$ and $S_r$ is a surface.  Then using Fujino and Mori's canonical bundle formula, see Lemma \ref{divisor},  we are reduced to looking at a big $\mathbb{Q}$-divisor $K_{S_r} +D$ on $S_r$.  The proof then naturally splits into two cases based on the geometry of $S_r$.  

When $\kappa(S_r) \geq 0$, the proof of Theorem \ref{dim3} reduces to a problem about the positivity of $K_{S_r} + D - \triangle$, where $\triangle$ is an effective $\mathbb{Q}$-divisor with simple normal crossing support.  From here, by passing to a minimal model of $S_r$, we may apply Reider's theorem, or a slight generalization of it, to conclude.

The other case we need to consider is when $\kappa(S_r) = - \infty$.  In this case, we are reduced to $S_r$ being a birationally ruled surface.  Using that the nef cone and the cone of curves on such a surface is relatively simple, we may directly prove that the integral divisor $(r-1)K_{S_r} + \lfloor r D \rfloor $ is nef.  Therefore, taking a sufficiently large multiple, which is effectively computable, we may once again apply Reider's theorem to conclude as in the previous case.  

We now give a brief outline of the paper.  In Section \ref{sectwo} we recall some background material along with our conventions.  In Section \ref{secthree} we prove Theorem \ref{dim3} when $\kappa(X)=2$.  Finally, in Section \ref{secfour} we prove the rest of the assertions in the introduction.

\noindent \emph{Acknowledgments:} It is a pleasure to thank my advisor M. Nakamaye for all of his help and encouragement.  I would like to thank E. Viehweg and D.-Q. Zhang for acknowledging my solution even though it came shortly after their stronger result.  I would also like to thank J. McKernan for suggesting this problem to me and G. Pacienza for giving a wonderful lecture series in Grenoble related to these ideas.  I would like to thank O. Fujino for pointing out two mistakes in an earlier version of this paper.  Finally, I would like to thank A. Buium and A. Saha for many useful discussions.

\section{Background}\label{sectwo}
We work over the field of complex numbers $\mathbb{C}$.  By a \emph{pair} $(X, \triangle)$, we mean a normal variety $X$ and an effective $\mathbb{Q}$-Weil divisor $\triangle$ such that $K_X + \triangle$ is $\mathbb{Q}$-Cartier.  We will say the pair $(X, \triangle)$ is \emph{Kawamata log terminal}, abbreviated klt, if there exists a log resolution $\mu : X' \to X$ of the pair $(X, \triangle)$ such that
$$
\mu^* (K_X + \triangle) = K_{X'} + E,
$$
where $E = \sum_{i} a_{i} E_i$ and each $E_i$ is a distinct irreducible divisor and every $a_i <1$.  Let $D = \sum_i a_i D_i$ be a $\mathbb{Q}$-divisor on a variety $X$, then the \emph{round-down} of $D$ is the integral divisor
$$
\lfloor D \rfloor = \sum_i \lfloor a_i \rfloor D_i,
$$
where if $x \in \mathbb{Q}$, then $\lfloor x \rfloor$ denotes the greatest integer less than or equal to $x$.  Similarly, we denote the \emph{round-up} by $\lceil x \rceil$.

By an \emph{algebraic fiber space}, we mean a surjective morphism $f: X \to S$ between  smooth projective varieties with connected geometric fibers.  Recall that any algebraic fiber space has the property that ${\phi_r}_* \mathscr{O}_X = \mathscr{O}_{S_r}$.  We now recall the notion of an Iitaka fibration:

\begin{theorem}
Let $X$ be a smooth projective variety with $\kappa(X) \geq 0$.  Then for all sufficiently large and divisible $r >0$ the rational maps
$$
\phi_r : X \dashrightarrow S_r \subset \mathbb{P} (H^0 (X, \mathscr{O}_X (r K_X)))
$$
are birationally equivalent to a fixed algebraic fiber space $\phi_{\infty} : X_{\infty} \to S_{\infty}$.  If $F_{\infty}$ denotes the generic fiber of $\phi_{\infty}$, then $\kappa(F_{\infty}) =0$.  Furthermore, $\kappa(X) = \dim S_{\infty}$.
\end{theorem}

\begin{proof}
We refer the reader to \cite[Defininition-Theorem 1.11]{Mori} or \cite[Theorem 2.1.33]{Lazarsfeld}.
\end{proof}

We call the algebraic fiber space $\phi_{\infty} : X_{\infty} \to S_{\infty}$ appearing in the previous theorem the \emph{Iitaka fibration} of $X$, which is well defined up to birational equivalence.  For other properties of algebraic fiber spaces along with basic properties of the Iitaka fibration, we refer the reader to Mori's survey paper \cite{Mori}.

We recall the following fundamental theorem of Reider, see \cite[Theorem 1]{Reider}, upon which our proof of Theorem \ref{dim3} heavily relies:

\begin{theorem}\label{basept}
Let $S$ be a smooth projective surface and $L$ a nef divisor on $X$.
\begin{enumerate}
\item If $L^2 \geq 5$ and the linear series $|K_S +L|$ has a base point $x \in S$, then there exists an effective divisor $C \subset S$ containing $x$ such that $L.C =0$ and $C^2 =-1$ or $L.C=1$ and $C^2 =0$.

\item If $L^2 \geq 10$ and the linear series $| K_S + L|$ fails to separate two points $x,y \in S$, then there exists an effective divisor $C \subset S$ containing $x$ and $y$ such that either $L.C =0$ and $C^2 =-1$ or $C^2 =-2$; or $L.C =1$ and $C^2 =0$ or $C^2 =-1$; or $L.C =2$ and $C^2 =0$.
\end{enumerate}
\end{theorem}

\begin{proof}
See \cite[Theorem 1]{Reider}.  For a nice introduction to the ideas surrounding this theorem see \cite[Theorem 2.1]{Lazarsfeld2}.
\end{proof}

Finally, we briefly recall some basic facts about ruled surfaces, see \cite[Section V.2]{Hartshorne} and \cite[Section 1.5.A]{Lazarsfeld}.  By a \emph{ruled surface}, we mean a smooth projective surface $S$ with a morphism $S \to C$ such that $C$ is a smooth projective curve of genus $g$ and every fiber of $S \to C$ is isomorphic to $\mathbb{P}^1$.  In this case, there exists a rank two vector bundle $\mathscr{E}$ on $C$ such that $S \cong \mathbb{P}(\mathscr{E})$.  Then after twisting $\mathscr{E}$ by a suitable line bundle we may suppose $\mathscr{E}$ has positive degree.  If $\mathscr{E}$ has a line bundle quotient of negative degree, then $S$ is called \emph{unstable}, otherwise we call $S$ \emph{stable}.  In either case we may write the N\'eron-Severi group of $S$ as $\textrm{Num}(X)_{\mathbb{R}} = \mathbb{R}[C_0] \oplus \mathbb{R}[f]$, where $[C_0]$ is the class of a section of $S \to C$ and $[f]$ is the class of a fiber of $S \to C$.

\section{When $\kappa(X)=2$}\label{secthree}
Let $X$ be a smooth projective threefold.  In this section we prove that Theorem \ref{dim3} is true in the case where $\kappa(X)=2$.  As was pointed out in the introduction the proof of this case breaks down into three cases based on the geometry of $S_r$, where $S_r = \phi_r (X)$ is the image of $X$ in $\mathbb{P} (H^0 (X, \mathscr{O}_X(rK_X)))$ for $r >>0$.

The main result of this section is:

\begin{theorem}\label{kap2}
Theorem \ref{dim3} is true when $\kappa(X)=2$.
\end{theorem}

We will first recall some background results, from \cite{MoriFuj}, used in the proof of Theorem \ref{kap2}.  Many of the methods we use are simialr to those found in \cite{MoriFuj}.  The main difference is that our fiber space is of relative dimension one instead of relative dimension two.  However, in our case we need to deal with the fact that our fibration $\phi_r : X \to S_r$ will have a higher dimensional base.  

After recalling some background results we will deal with the easier of the two cases, namely, the case where $\kappa(S_r) \geq 0$, and afterward we will deal with the case when $\kappa(S_r) =- \infty$.

\begin{lemma}\label{ell}
It suffices to prove Theorem \ref{kap2}, when $\phi_r : X \to S_r$ is an elliptic fiber space.
\end{lemma}

\begin{proof}
Let $\phi_{\infty} : X _{\infty} \to S_{\infty}$ be the Iitaka fibration.  By \cite[Remark 1.11.1(b)]{Mori} the generic fiber of $\phi_r : X \dashrightarrow S_r$ will be a non-singular elliptic curve.  Furthermore, $X$ will be birational to $X_{\infty}$ and $S_r$ will be birational to $S_{\infty}$.  By making a birational modification of $X$ we may assume that $\phi_r : X \to S_r$ is a proper morphism and that $X$ and $S_r$ are both non-singular.  Finally, by \cite[Proposition 1.4]{Mori} $\phi_r : X \to S_r$ will be an elliptic fiber space.
\end{proof}

Throughout the rest of this section $\phi_r : X \to S_r$ will be an elliptic fiber space, unless we say otherwise, with $X$ and $S_r$ non-singular with generic fiber $F$, which is necessarily an elliptic curve.  

We will need the following result of Fujino and Mori, see \cite[Proposition 2.2]{MoriFuj}, whose proof we recall for the convenience of the reader.

\begin{lemma}\label{divisor}
Let $\phi_r : X \to S_r$ be an elliptic fiber space.  Then up to linear equivalence there exists a $\mathbb{Q}$-divisor $D$ on $S_r$ such that
\begin{equation}\label{blahblah}
K_X = \phi_r ^* ( K_{S_r} +D) +B,
\end{equation}
where $B$ is a $\mathbb{Q}$-divisor on $X$.  Furthermore, if we write $B = B^+ - B^-$ as the difference of two effective $\mathbb{Q}$-divisors without common components, then we have ${\phi_r}_* \mathscr{O}_X ( \lfloor mB^+ \rfloor) = \mathscr{O}_{S_r}$ for every $m >0$ and  $\Supp ( \phi_r (B^{-}))$ has codimension two in $S_r$.
\end{lemma}

\begin{proof}
By \cite[Theorem 2.6(i)]{Mori} there exists a $c \in \mathbb{N}$ such that
$$
({\phi_r}_* \mathscr{O}_X (mc K_{X/ S_r}))^{**} \cong (({\phi_r}_* \mathscr{O}_X (c K_{X/S_r}))^{**})^{\otimes m}
$$
for all $m \in \mathbb{N}$, where ${}^{**}$ denotes the double dual.  Indeed, in our restricted setting it suffices to take $c$ to be even.  Since $S_r$ is smooth and ${\phi_r}_* \mathscr{O}_X (c K_{X/S_r})^{**}$ is an invertible sheaf there exists a $\mathbb{Q}$-divisor $D$ on $S_r$ such that
$$
\mathscr{O}_{S_r} (cD) \cong {\phi_r}_* \mathscr{O}_X (c K_{X/S_r})^{**}.
$$
By construction it is unique up to linear equivalence.  Since taking the double dual of ${\phi_r}_* \mathscr{O}_X (bc K_{X/S_r})$ does not effect the divisor $D$ associated to this sheaf we see that
$$
\phi_r ^* {\phi_r}_* \mathscr{O}_X (c K_{X /S_r}) \hookrightarrow \mathscr{O}_X (c K_{X/S_r}) |_U,
$$
where $U$ is $X$ minus some codimension two algebraic set of $X$.  Therefore, we may take $K_{X /S_r} - \phi_r ^* D = B$ and it is easily checked that $B$ satisfies the stated hypothesis.
\end{proof}

The proof of Lemma \ref{divisor} actually holds in a more general setting.  Namely, when the generic fiber $F$ of $\phi_r : X \to S_r$ has dimension greater than one and $\kappa(F)=0$.  In this more general setting we must replace (\ref{blahblah}) by
$$
bK_X = \phi_r ^* ( bK_{S_r} +D) +B,
$$
where $b = \min \{ b \in \mathbb{N} : H^0 (F, \mathscr{O}_F(bK_F)) \neq 0 \}$.  We should also mention that the $\mathbb{Q}$-divisor $D$ appearing the previous lemma is the unique $\mathbb{Q}$-divisor, up to numerical equivalence, which statisfies the $\mathscr{O}_{S_r}$-algebra isomorphism
$$
\oplus_{m \in \mathbb{N}} \mathscr{O}_{S_r} (\lfloor mD \rfloor) \cong \oplus_{m \in \mathbb{N}} \mathscr{O}_{S_r} (\phi_{r*} m b K_{X/S_r})^{**}.
$$

By a fundamental result of Fujino and Mori \cite[Proposition 2.8]{MoriFuj} it is possible to give a rather explicit description of the divisor $D$ in the previous lemma.  Fujita has also proved a similar statement \cite[Theorem 2.15]{Fujita}.

\begin{lemma}\label{bignef}
We may write the divisor $D$ appearing in the previous lemma as 
$$
D = L_r + \sum_{i} s_{i} E_i,
$$
where $L_r$ is nef $\mathbb{Q}$-divisor and $\triangle = \sum_i s_{i} E_i$ is  an effective $\mathbb{Q}$-divisor with $ \lfloor \triangle \rfloor =0$.
\end{lemma}

\begin{proof}
We may write $D= L_r + \triangle$ by \cite[Proposition 2.8]{MoriFuj}, where $\triangle =\sum_{i} s_i E_i$ is effective.  Furthermore, the terms $s_i$ are of the form
$$
s_{i} = \frac{N u_i - v_i}{N u_i},
$$
where $u_i, v_i \in \mathbb{N}$ with $0 < v_i \leq N$ and $N \in \mathbb{N}$ is the smallest integer such that $N \phi_r ^* L_r$ is an integral divisor.  The above immediately implies $ \lfloor  \triangle  \rfloor =0$.  The fact that $L_r$ is nef follows from the fact that $L_r$ is the pull-back of $\mathscr{O}_{\mathbb{P}^1} (1)$ by a suitable projective morphism \cite[Theorem 2.15]{Fujita}.  By Lemma \ref{divisor}, $D$ is unique up to linear equivalence.  Thus, it must coincide with the canonical bundle formula of Fujino and Mori, see \cite[Theorem 2.15]{Fujita}.  Then $NL_r = J^* \mathscr{O}_{\mathbb{P}^1} (1)$, where $N$ is the smallest natural number for which $NL_r$ is an integral divisor and $J : S_r \to \mathbb{P}^1$ is the extension of the $j$-invariant of the fibers of $\phi_r$ on a dense open set $U \subset S_r$ to all of $S_r$, after possibly making a birational modification to $S_r$.  
\end{proof}

The following result of Fujino and Mori, see \cite[Theorem 3.1]{MoriFuj}, gives a criterion for when the divisor $L_r$ will be an integral divisor.  They also point out that an affirmative answer for Conjecture \ref{basic}, when $r$ is sufficiently divisible, could be obtained in the case when $\dim X =n$ and $\kappa(X)=1$, if one could effectively bound the Betti number $B_{n-1}$ of the generic fiber of $\phi_{\infty} : X_{\infty} \to S_{\infty}$. 

\begin{lemma}\label{betti}
Let $\phi_r : X \to S_r$ be an elliptic fiber space then $N L_r$ is an integral divisor whenever $N \in \mathbb{Z}$ is divisble by $12$.
\end{lemma}

\begin{proof}
This is a special case of \cite[Theorem 3.1]{MoriFuj}.  The proof shows that we may take $N (B_1) = \textrm{lcm} \{ x \in \mathbb{N} : \varphi (x) \leq B_1 \}$ where $\varphi$ is Euler's function and $B_1$ is the first Betti number of the generic fiber.  Hence, the divisor $NL_r$ will be integral whenever $N$ is divisible by $12$.
\end{proof} 

The reason that our proof of Theorem \ref{dim3}, in the case when $\kappa(X)=2$, only holds for $r$ sufficiently divisible is because we need $L_r$ to be an integral divisor when we take cohomology since it is not clear how the divisor $L_r$ behaves when we take the round-down.  Throughout the rest of the section we make the assumption that $r$ is sufficiently divisible, which by the previous lemma means that we require $r$ to be divisible by $12$.

The next lemma shows that the singularities of the divisor $\triangle$ are suitably mild.  This fact will play a crucial role in the later proofs since one can apply various minimal model programs to surfaces in a fairly controlled way.

\begin{lemma}\label{klt}
Let $\phi_r : X \to S_r$ be an elliptic fiber space and $\triangle$ the $\mathbb{Q}$-divisor appearing in Lemma \ref{bignef}.  Then the pair $(S, \triangle)$ is klt.
\end{lemma}

\begin{proof}
By Lemma \ref{bignef}, the divisor $\triangle$ is effective with $ \lfloor \triangle  \rfloor=0$.  Therefore, by \cite[Remark 5.15.9(i)]{Mori}, the divisor $\triangle$ is a negligible divisor, recall this implies $\Supp (\triangle)$ is a simple normal crossing divisor.  By \cite[Corollary 2.3.1(3)]{KollarMori}, it follows that $\textrm{discrep}(S_r, \triangle) >-1$ since $S_r$ is smooth and hence the pair $(S_r, \triangle)$ is klt. 
\end{proof}

\begin{lemma}
Let $\phi_r : X \to S_r$ and $L_r$ be as in Lemma \ref{bignef}.  To prove Theorem \ref{dim3} is true when $\dim X =3$, and the Iitaka fibration stabilizes to 
$$
\phi_{\infty} : X _{\infty} \dashrightarrow S_{\infty},
$$
where $S_{\infty}$ is a surface with $\kappa(S_{\infty}) \geq 0$, it suffices to prove there exists a universal constant $r_3 $ such that the natural map
$$
\psi_r : S_r \to \mathbb{P}(H^0 (S_r , \mathscr{O}_{S_r} (r K_{S_r} + rL_r))),
$$
is a birational morphism onto its image whenever $r \geq r_3 $ and sufficiently divisible.
\end{lemma}

\begin{proof}
By Lemma \ref{ell} we may assume that $\phi_r : X \to S_r$ is an elliptic fiber space.  By Lemma \ref{divisor} we may write
\begin{equation}\label{form}
K_X = \phi_r ^* ( K_{S_r} +D) +B,
\end{equation}
where $D=L_r + \triangle$, with $L_r$ nef.  Since ${\phi_r}_* \mathscr{O}_X ( \lfloor rB^+ \rfloor) = \mathscr{O}_{S_r}$  and $\Supp(B^-)$ has codimension two in $S_r$ we see by (\ref{form}) that
$$
H^0 (X, \mathscr{O}_X (r K_X)) = H^0 (X, \mathscr{O}_X (\lfloor \phi_r ^* (r K_{S_r} + rD) \rfloor ))
$$
for every $r >0$, by \cite[Theorem 4.5(iii)]{MoriFuj} .  Therefore, it follows that
\begin{equation}\label{final}
H^0 (X, \mathscr{O}_X (r K_X)) = H^0 (S_r, \mathscr{O}_{S_r} ( \lfloor r K_{S_r} + rD \rfloor )).
\end{equation}
This implies that the $\mathbb{Q}$-divisor $K_{S_r} +D$ is big.  For $r$ divisible by $12$ the divisor $rL_r$ is integral by Lemma \ref{betti}.  Therefore, we have
$$
\lfloor r K_{S_r} + rL_r + r \triangle  \rfloor = r K_{S_r} + rL_r + \lfloor r \triangle \rfloor.
$$
Since $\lfloor r \triangle \rfloor$ is effective, we get an injective map
$$
H^0 (S_r , \mathscr{O}_{S_r} (rK_{S_r} + r L_r)) \hookrightarrow H^0 (S_r , \mathscr{O}_{S_r} (\lfloor r K_{S_r} + rL_r + r \triangle  \rfloor )).
$$
If $H^0 (S_r, \mathscr{O}_{S_r} (r K_{S_r} + r L_r))$ induces a birational morphism
$$
\psi_r : S_r \to \mathbb{P} (H^0 (S_r, \mathscr{O}_{S_r} (r K_{S_r} +r L_r)))
$$
for all $r \geq r_3 $ and sufficiently divisible, where $r_3 $ is the universal constant appearing in the statement of the lemma, then by the bijection (\ref{final}) we see that
$$
\phi_r : X \to S_r \subset \mathbb{P}(H^0 (X, \mathscr{O}_{X} (r K_X)))
$$
will be birational to the Iitaka fibration and have connected fibers.
\end{proof}

\begin{lemma}\label{thingy23}
Theorem \ref{dim3} is true when $\dim X =3$, and the Iitaka fibration stabilizes to 
$$
\phi_{\infty} : X _{\infty} \dashrightarrow S_{\infty},
$$
where $S_{\infty}$ is a surface with $\kappa(S_{\infty}) \geq 0$ and $K_{S_r} + L_r$ is big.
\end{lemma}

\begin{proof}
With the notation as in the previous lemma, it suffices to show that there exists a constant $r_3$ such that $H^0 (S_r, \mathscr{O}_{S_r} (r K_{S_r} + r L_r ))$ induces a birational morphism
$$
\psi_r : S_r \to \mathbb{P} (H^0 (S_r, \mathscr{O}_{S_r} (r K_{S_r} + r L_r ))),
$$
whenever $r \geq r_3$ and sufficiently divisible.  Let $S_r '$ be the minimal model of $S_r$ and $f : S_r \to S_r '$ the birational morphism.  Then $S_r '$ is smooth as $S_r$ is smooth and $\dim S_r =2$.  Since $f$ is birational with connected fibers, we have
$$
H^0 (S_r , \mathscr{O}_{S_r} ( r K_{S_r} + r L_r )) = H^0 (S_r ' , \mathscr{O}_{S_r '} (r K_{S_r ' } + f_* (rL_r))).
$$
Since $S_r '$ is minimal, we see that $(r-1) K_{S_r '} + f_* (r L_r)$ is big and nef.  By possibly replacing $S_r$ by $S_r '$ we may assume that $S_r$ is a smooth minimal surface and $rK_{S_r} + r L_r$ is big and nef, where the bigness of our divisor follows from our hypothesis.  For $r \geq 48$ and sufficiently divisible, we have
$$
 ((r-1) K_{S_r} + r L_r)^2 \geq 10.
 $$
To see this notice that all the intersection numbers appearing in the intersection are non-negative.  By Theorem \ref{basept}, if the linear series $|r K_{S_r} + r L_r|$ fails to separate points, then there exists an effective divisor $C \subset S_r$ such that $((r-1) K_{S_r} + r L_r).C =0$ since our hypothesis on $r$ eliminates the other possible cases appearing in Theorem \ref{basept}.  Since $L_r$ is nef, we must have $K_{S_r} .C =0$.  By adjunction and the fact that $S_r$ is minimal, we see that the morphism
$$
\psi_r : S_r \to  \mathbb{P} (H^0 (S_r, \mathscr{O}_{S_r} (r K_{S_r} + rL_r)))
$$
defines an embedding outside of $(-2)$-curves, again by Theorem \ref{basept}.  Applying the Stein factorization to $\psi_r$ shows that $\psi_r$ has connected fibers and hence is a birational morphism onto its image.
\end{proof}

\begin{proposition}\label{one}
Theorem \ref{dim3} is true when $\dim X =3$, and the Iitaka fibration stabilizes to 
$$
\phi_{\infty} : X _{\infty} \dashrightarrow S_{\infty},
$$
where $S_{\infty}$ is a surface with $\kappa(S_{\infty}) \geq 0$.
\end{proposition}

\begin{proof}
We may assume that $K_{S_r} + L_r$ is nef, by a similar construction, as in the proof of the previous lemma.  Furthermore, we may assume that $K_{S_r}+L_r$ is not big.  If we run a log-MMP on the pair $(S_r, \triangle)$ then by \cite[Theorem 2.3.6]{KollarKovacs} the end result will be a pair $(S_r ', \triangle ')$ with $S_r '$ birational to $S_r$ and $K_{S_r ' } + \triangle '$ big and nef, as $K_{S_r} + \triangle$ is big.  Since $\lfloor \triangle ' \rfloor =0$ we may apply \cite[Theorem 4.5]{KollarMori} and \cite[Theorem 4.7]{KollarMori} to conclude that $S_r '$ is in fact smooth, since $\triangle$ is a simple normal crossings divisor and each component is a smooth rational curve.  That is, the pair $(S_r ' , \triangle ')$ is in fact canonical and by the explicit description of $\triangle$, given in Lemma \ref{klt}, we conclude that $S_r '$ is smooth.  Then replacing $(S_r, \triangle)$ with the pair $(S_r ' , \triangle ')$ we obtain an algebraic fiber space $ \phi_r : X \to S_r$, where $K_{S_r} + L_r + \triangle$ is big and nef.  For $r$ divisible by $12$ we see that the integral divisor
$$
(r-1)K_{S_r} + r L_r + \lceil r \triangle \rceil
$$
is big.  If $\lceil r \triangle \rceil$ is nef, then we may apply Theorem \ref{basept} to conclude, as we did in the proof of Lemma \ref{thingy23}.  On the other hand if 
$$
(r-1) K_{S_r} + r L_r + \lceil r \triangle \rceil
$$
is not nef and hence $\lceil r \triangle \rceil$ is not nef, then we must use a generalization of Theorem \ref{basept}, see \cite[Proposition 4]{Masek}, to deal with the non-nefness of $\lceil r \triangle \rceil$.  Indeed, let $r \geq 36$ and divisible by $12$ then if $|r K_{S_r} + r L_r + \lceil r \triangle \rceil |$ fails to seperate two points, there exists an irreducible curve $C \subset S_r$ with $K_{S_r} . C =0$ and $ \lceil r \triangle \rceil . C  \leq 0$.  It then follows that $C \subset \textrm{Supp}(\triangle)$ and hence rational by \cite[Theorem 4.7]{KollarMori}.  Therefore, by adjunction, we conclude that $C^2 = -2$.  Hence, the map
$$
\psi_r : S_r \to \mathbb{P}(H^0 (S_r, \mathscr{O}_{S_r} ( rK_{S_r} + r L_r + \lceil r \triangle \rceil )))
$$
is an embedding outside of $(-2)$-curves and our claim follows.
\end{proof}

We now deal with the case where the Iitaka fibration stabilizes to $\phi_{\infty} : X_{\infty} \to S_{\infty}$ where $\kappa(S_{\infty}) = - \infty$.  First we need the following lemma.

\begin{lemma}
Suppose that the Iitaka fibration stabilizes to 
$$
\phi_{\infty} : X_{\infty} \dashrightarrow S_{\infty}
$$
where $\kappa(S_{\infty}) = - \infty$.  It suffices to prove Theorem \ref{dim3} when $S_r$ is a ruled surface.
\end{lemma}

\begin{proof}
As in the proof of Lemma \ref{one} we may assume $\phi_r : X \to S_r$ is an elliptic fiber space and 
\begin{equation}\label{thingy}
H^0 (X, \mathscr{O}_X (r K_X)) = H^0 (S_r, \mathscr{O}_{S_r} ( \lfloor r K_{S_r} + rD \rfloor )),
\end{equation}
where $D$ is as in Lemma \ref{bignef}.  Consider the pair $(S_r, \triangle)$.  Running a log-MMP program on this pair, by \cite[Theorem 2.3.6]{KollarKovacs}, we obtain a pair $(S_r ', \triangle ')$, where $S_r '$ is a birationally  ruled surface or $(S_r ' , \triangle ')$ is a log Del Pezzo surface.  In either case $S_r '$ will be birational to a smooth ruled surface $S_r ''$.  Let $\psi : S_r \to S_r ''$ be the birational morphism to the smooth ruled surface $S_r ''$.  After possibly replacing $S_r$ with $S_r ''$ and $\phi_r$ with $\psi \circ \phi_r$ we may assume that $\phi_r : X \to S_r$ is an elliptic fiber space and $S_r$ is a ruled surface.  By abuse of notation we replace $L_r$ with $\psi_{*} L_r$, which is a nef $\mathbb{Q}$-divisor on $S_r$. 
\end{proof}

\begin{proposition}\label{rulesurf}
Theorem \ref{dim3} is true when $\dim X =3$, and the Iitaka fibration stabilizes to 
$$
\phi_{\infty} : X _{\infty} \dashrightarrow S_{\infty},
$$
where $S_{\infty}$ is rational.
\end{proposition}

\begin{proof}
By the previous lemma, we may assume $\phi_r : X \to S_r$ is an elliptic fiber space with $S_r \to C$ a ruled surface.  We denote the genus of $C$ by $g$.  Then we may write $S_r = \mathbb{P} (\mathscr{E})$ for a suitable rank two vector bundle $\mathscr{E}$ on $C$ where, after twisting by a suitably positive line bundle, we may assume $\deg \mathscr{E}=-e \geq 0$.  By \cite[Corollary V.2.11]{Hartshorne} we may write $K_{S_r}$, up to numerical equivalence, as
\begin{equation}\label{can}
K_{S_r } = -2C_0 +(2g -2 -e)f,
\end{equation}
where $C_0$ is the divisor corresponding to a section of $S_r  \to C$ with $C_0 ^2 =-e \geq 0$ and $f$ denotes the class of a fiber of $S_r$.  We claim that the integral divisor $  (r-1)K_{S_r} + rL_r + \lfloor r \triangle \rfloor$ is nef.  To check this, we verify the claim in two cases based on whether the ruled surface $S_r \to C$ is stable or unstable.  We refer the reader to \cite[Section 1.5.A]{Lazarsfeld} for a description of the various cones that can occur in each case.

First suppose that $S_r \to C$ is stable.  In this case, the closed cone of curves $\overline{\textrm{NE}}(S_r)$ is equal to the nef cone $\textrm{Nef}(S_r)$.  Since $L_r$ is nef and $\textrm{Num}(S_r)_{\mathbb{R}}$ is generated by $C_0$ and $f$, we must have
$$
rL_r =  aC_0 +bf,
$$
where $a \geq 0$ and $aC_0 ^2 + b \geq 0$.  To see this recall that the intersection numbers are $C_0 ^2 =-e$, $f^2 =0$ and $C_0 . f =1$.  By using (\ref{can}), we have
$$
(r-1) K_{S_r} + rL_r = (a-2(r-1))C_0 + (2g-2 -e)(r-1)f + bf.
$$
Since $\lfloor r \triangle \rfloor$ is effective, we may write $\lfloor r\triangle \rfloor= \lfloor \alpha   \rfloor C_0+  \lfloor \beta \rfloor f$ with $ \alpha, \beta \geq 0$.  Therefore,
\begin{equation}\label{firstint}
( (r-1) K_{S_r} + r L_r +  \lfloor r \triangle \rfloor ).C_0 = (a-2(r-1) + \lfloor \alpha \rfloor ) C_0 ^2 + b +(2g-2 -e)(r-1) + \lfloor \beta \rfloor 
\end{equation}
and
\begin{equation}\label{secondint}
(K_{S_r} + rL_r + \lfloor r\triangle \rfloor ).f = a -2(r-1) + \lfloor \alpha \rfloor .
\end{equation}
Since $K_{S_r} + L_r + \triangle $ is pseudoeffective, so is $ (r-1) K_{S_r} + r L_r + \lfloor r \triangle \rfloor$ and so it has non-negative intersection with all nef divisors since $S_r$ is a surface, see \cite[Remark 2.2.27]{Lazarsfeld}.  Therefore, both (\ref{firstint}) and (\ref{secondint}) will be non-negative, implying that $ (r-1) K_{S_r} + rL_r + \lfloor r\triangle \rfloor$ is nef since in the case when $S_r$ is stable the nef cone and the closure of the cone of effective curves coincide.

We now suppose that  the ruled surface $S_r \to C$ is unstable.  In this case the the two cones $\overline{\textrm{NE}}(S_r)$ and $\textrm{Nef}(S_r)$ differ.  However, both are two-dimensional.  Let $-\delta f+C_0$ and $f$ denote the two rays bounding $\textrm{Nef}(S_r)$, where $\delta <0$.  Then we may write $rL_r = aC_0 +bf$ with $a \geq 0$ and $b \geq 0$.  To see this, use the fact that $rL_r$ is nef so that $r L_r . f =a \geq 0$ and $rL_r .(\delta f + C_0) = a C_0 ^2 + \delta a + b \geq 0$.  Similarly, since $\triangle$ is an effective $\mathbb{Q}$-divisor we have $\lfloor r\triangle \rfloor = \lfloor \alpha \rfloor C_0 + \lfloor \beta \rfloor f$ with $ \lfloor r \triangle \rfloor .f = \lfloor \alpha \rfloor \geq 0$ and $ \lfloor r \triangle \rfloor . (-\delta f + C_0) = -\delta  \lfloor \alpha \rfloor  + aC_0 ^2+ \lfloor \beta \rfloor \geq 0$.  Thus, up to numerical equivalence, we may write
$$
(r-1) K_{S_r} +r L_r +  \lfloor r \triangle \rfloor = (a+ \lfloor \alpha \rfloor -2(r-1))C_0 + (b + \lfloor \beta \rfloor + (2g-2 -e)(r-1) )f.
$$
Now suppose that $(r-1) K_{S_r} + r L_r + \lfloor r \triangle \rfloor$ is not nef.  Then there exists an effective curve, which after scaling we may write as $C' = \delta ' f + C_0$ where $\delta ' \geq \delta$, with 
\begin{equation}\label{ifzero}
(r-1)K_{S_r} +r L_r + \lfloor r \triangle \rfloor . C' =(a + \lfloor \alpha \rfloor -2(r-1) )C_0 ^2 + \delta ' (a+ \lfloor \alpha \rfloor  -2(r-1)) + b +  \lfloor \beta \rfloor + (2g-2 -e)(r-1) < 0.
\end{equation}
Since $(r-1) K_{S_r} +  rL_r + \lfloor r \triangle \rfloor$ is pseudoeffective, $C '$ is not a nef divisor.  By rearranging the left side of (\ref{ifzero}), we have 
\begin{equation}\label{thingy}
\frac{(a + \lfloor \alpha \rfloor -2(r-1))C_0 ^2+ b +  \lfloor \beta \rfloor + (2g-2-e)(r-1)} {a+  \lfloor \alpha \rfloor  -2 (r-1) } <- \delta ' \leq - \delta
\end{equation}
since $-e >0$.  Notice that if the denominator is zero, then the divisor in question is already nef.  Therefore, we may assume,
$$
( (r-1) K_{S_r} +rL_r + \lfloor r \triangle \rfloor ).f = a + \lfloor \alpha \rfloor -2(r-1) > 0.
$$
Since $\delta f + C_0$ and $f$ bound $\overline{\textrm{NE}}(S_r)$, using (\ref{thingy}), we see that
$$
 0 \leq \frac{b +  \lfloor \beta \rfloor +(2g-2-e)(r-1)} {a+  \lfloor \alpha \rfloor  -2 (r-1) } \leq - \delta
$$
since $(r-1) K_{S_r} + r L_r + \lfloor r \triangle \rfloor$ is psuedoeffective.  But, this implies that $  (r-1) K_{S_r} + rL_r+ \lfloor r\triangle \rfloor$ lies in $\textrm{Nef}(S_r)$, which is a contradicition.  Therefore, the integral divisor $(r-1) K_{S_r} + rL_r + \lfloor  r \triangle \rfloor$ is nef.

To conclude we may now apply Theorem \ref{basept} to the integral divisor $ (r-1) K_{S_r} + rL_r + \lfloor r\triangle \rfloor$ with $r \geq 48$ and sufficiently divisible.  By our hypothesis on $r$, it is easily verified that $ | rK_{S_r} + r L_r + \lfloor r \triangle \rfloor |$ is base point free and that the linear series separates points of $S_r$ outside of $(-2)$-curves.  Therefore, we may conclude that the morphism
$$
 \psi_r : S_r \to \mathbb{P}(H^0 (S_r, \mathscr{O}_{S_r} (r K_{S_r} + r L_r +\lfloor r \triangle \rfloor)))
$$
is birational onto its image with connected fibers, as in the proof of Proposition \ref{one}.  Hence, $\phi_r : X \to S_r$ is birational to the Iitaka fibration whenever $r \geq 48$ and sufficiently divisible.
\end{proof}

\begin{proof}[Proof of Theorem \ref{kap2}]
By Lemma \ref{ell}, it suffices to prove the proposition when $\phi_r : X \to S_r$ is an elliptic fiber space.  If the Iitaka fibration stabilizes with $\kappa(S_{\infty}) \geq 0$, then we may apply Proposition \ref{one} to conclude that $\phi_r : X \to S_r$ is birational to the Iitaka fibration as soon as $r \geq 48$ and sufficiently divisible, where by sufficiently divisible we mean divisible by $12$.  If $\kappa(S_{\infty}) = -\infty$, then we may apply Proposition \ref{rulesurf} and $\phi_r : X \to S_r$ will be birational to the Iitaka fibration as soon as $r \geq 48$ and divisible by $12$.  Hence, the claim follows.
\end{proof}

\section{Proof of Theorem \ref{dim3}}\label{secfour}

Based on the results of Section \ref{secthree} along with the work of Hacon, Fujino, McKernan, and Mori, the proof of Theorem \ref{dim3} is now trivial.

\begin{proof}[Proof of Theorem \ref{dim3}]
The case where $\kappa(X)=0$ follows from \cite[10.4.3]{Mori}.  The case when $\kappa(X) =1$ was dealt with by Fujino and Mori \cite[Corollary 6.2]{MoriFuj}.  The case where $\kappa(X) =2$ is the content of Theorem \ref{kap2}.  Finally, the case where $\kappa(X)=3$ was dealt with by Hacon and McKernan \cite[Theorem 1.1]{McKernan}.  
\end{proof}

\begin{proof}[Proof of Corollary \ref{kap2bound}]
This follows directly from the proofs of Proposition \ref{one} and Proposition \ref{rulesurf}.
\end{proof}

\begin{proof}[Proof of Corollary \ref{explicit}]
In the case when $\kappa(X)=3$, by \cite[Theorem 1.1]{ChenChen}, it suffices to take $r \geq 77$.  In the case where $\kappa(X)=2$, we have seen by Corollary \ref{kap2bound} that it suffices to take $r \geq 48$ and sufficiently divisble.  In the case where $\kappa(X)=1$, use \cite[Proposition 6.3]{MoriFuj}, which gives an effectively computable number.  Finally, for the case when $\kappa(X)=0$, see \cite[10.4.3]{Mori}.
\end{proof}


\end{document}